\documentclass[12pt]{article}

\usepackage{amsfonts}

\newtheorem{thm}{Theorem}[section]
\newtheorem{lemma}[thm]{Lemma}
\newtheorem{cor}[thm]{Corollary}

\newtheorem{prop}[thm]{Proposition}

\newenvironment{remark}{\par\medskip\noindent{\bf Remark.\ }}{\par\smallskip}
\newcommand{\proof
}{\par\medskip\noindent {\bf Proof.\ \ }}

\newcommand{\be}{\begin{equation}}
\newcommand{\ee}{\end{equation}}
\newcommand{\openbox}{\leavevmode
  \hbox to8pt{\hfil\vrule\vbox to6pt{\hrule width6pt\vfil\hrule}\vrule}}

\newcommand{\qed}{\hbox to5pt{ } \hfill \openbox\bigskip\medskip}

\newcommand{\cF}{\mbox{$\cal F$}}

\newcommand{\Z}{\mathbb Z}
\newcommand{\Q}{\mathbb Q}

\title{Betti numbers of Stanley--Reisner rings with pure resolutions}
\author{G\'abor Heged\"{u}s
\\{\normalsize Johann Radon Institute for Computational and Applied Mathematics}
}

\begin{document}


\maketitle

\begin{abstract}
Let $\Delta$ be simplicial complex and let $k[\Delta]$ denote the 
Stanley--Reisner ring corresponding to $\Delta$. 
Suppose that $k[\Delta]$ has a pure free resolution. 
Then we describe the Betti numbers and the Hilbert--Samuel multiplicity 
of $k[\Delta]$ in terms of the $h$--vector 
of $\Delta$. As an application, we derive a linear equation system 
and some inequalities for the components of 
the $h$--vector
of the clique complex of an arbitrary chordal graph.
As an other application, we derive a linear equation system 
and some inequalities for the components of 
the $h$--vector of Cohen--Macaulay simplicial complexes.
\end{abstract}
\medskip


\section{Introduction}
\noindent

Let $k$ denote an arbitrary field.
Let $R$ be the graded ring $k[x_1,\ldots ,x_n]$. 
The vector space $R_s=k[x_1,\ldots ,x_n]_s$ consists of 
the homogeneous polynomials of total degree $s$, together with $0$.

Let $G$ be a simple graph. 
We may think of an edge $E$ of a graph as a squarefree monomial
 $x^{E}:=\prod_{j\in E} x_j$ in $R$.

The {\em edge ideal} $I(G)$ is the ideal $\langle x^{E}:~ E\in E(G)\rangle$, 
which is generated by the edges of $G$.

The edge ideal was first introduced by R. Villareal in \cite{V}.
Later edge ideals have been studied very widely, see for instance 
\cite{E, E2, F, F2, Fr, HT, HT2, MRV, V, W, Z}.

In \cite{Fr} R. Fr\"oberg characterized the graphs $G$ such that $G$ 
has a linear free resolution. He proved:

\begin{thm} \label{Froberg_theorem}
Let $G$ be  a simple graph on $n$ vertices. Then $R/I(G)$ has linear 
free resolution precisely when $\overline{G}$, the complementary graph of $G$ 
is chordal. 
\end{thm}

In \cite{E2} E. Emtander generalized Theorem \ref{Froberg_theorem} for 
generalized chordal hypergraphs.

Let $\Delta$ be a simplicial complex. A facet $F$ is called a {\em leaf}, if 
either $F$ is the only facet of $\Delta$, or there exists an other facet
$G$, $G\neq F$ such that $H\cap F\subset G\cap F$ for each facet $H$ with 
$H\neq F$. 
A facet $G$ with this property is called a {\em branch} of $F$. 

In \cite{HHMTZ} J. Herzog, T. Hibi, S.  Murai, N. V. Trung and X. Zheng
 call the simplicial complex $\Delta$ a {\em quasi--forest} if
there exists a labeling $F_1,\ldots ,F_m$ of the facets such that 
for all $i$ the facet $F_i$ is a leaf of the subcomplex
 $\langle F_1,\ldots ,F_i\rangle$. We call such a labeling 
a {\em leaf order}. 

A graph is called {\em chordal} if each cycle of length $>3$ has a chord.
 
J. Herzog, T. Hibi, S.  Murai, N. V. Trung,  X. Zheng proved the following equivalent 
characterization of quasi--forests in \cite{HHMTZ} Theorem 1.1:
\begin{thm} \label{quasi_forest}
Given a finite sequence $(f_{-1},\ldots ,f_{d-1})$ of integers with each $f_i>0$,
the following conditions are equivalent:\\
(i) There is a quasi--forest $\Delta$ of dimension $d-1$ with
$f(\Delta):=(f_{-1},\ldots ,f_{d-1})$. \\
(ii) The sequence $(c_1,\ldots ,c_d)$, defined by the formula
$$
\sum_{i=0}^d f_{i-1}(x-1)^i=\sum_{i=0}^d c_ix^i,
$$ 
where $f_{-1}=1$, satisfies $\sum_{i=0}^d c_i>0$ for each $1\leq k\leq d$.\\
(iii) The sequence $(b_1,\ldots ,b_d)$  defined by the formula
$$
\sum_{i=0}^d f_{i-1}(x-1)^i=\sum_{i=0}^d b_ix^{i-1},
$$ 
is positive, i.e., $b_i>0$ for each $1\leq i\leq d$.
\end{thm}

Let $G$ be a finite graph on $[n]$. A {\em clique} of $G$ is a subset 
$F$ of $[n]$ such that $\{i,j\}\in E(G)$ for all $i,j\in F$ with $i\neq j$.

We write $\Delta(G)$ for the simplicial complex on $[n]$
whose faces are the cliques of $G$. This $\Delta(G)$ is the {\em clique complex} 
of the graph $G$.

We recall here for the famous  Dirac's Theorem (see \cite{D}).
\begin{thm} \label{Dirac} (Dirac)
A finite graph $G$ on $[n]$ is a chordal graph iff $G$ is the $1$--skeleton 
of a quasi--forest
\end{thm}

Let $\Delta$ be a simplicial complex. Our main results are 
some explicit formulas for the Betti numbers 
of the Stanley--Reisner ring $k[\Delta]$
such that $k[\Delta]$ has a pure free resolution in terms of the $h$--vector
of $\Delta$. 
 
As an application we give some linear equations 
for the components of the $h$--vector
of the clique complex of an arbitrary chordal graph.

In Section 2 we collected  some basic results about simplicial complices,
 free resolutions, Hilbert fuctions and Hilbert series.
 We present our main results in Section 3. We give some applications 
for Cohen--Macaulay squarefree monomial ideals and chordal graphs 
in Section 4.


\section{Preliminaries}
\subsection{Free resolutions}

Recall that for every finitely generated graded module $M$ over $R$ 
we can associate to $M$ 
a {\em minimal graded free resolution} 
$$ 
0\longrightarrow \bigoplus_{i=1}^{\beta_p} R(-d_{p,i}) \longrightarrow 
\bigoplus_{i=1}^{\beta_{p-1}} R(-d_{p-1,i})\longrightarrow \ldots \longrightarrow 
\bigoplus_{i=1}^{\beta_{0}} R(-d_{0,i}) \longrightarrow M \longrightarrow
0, 
$$
where $p\leq n$ and $R(-j)$ is the free $R$-module obtained 
by shifting the degrees of $R$ by $j$. 

Here the natural number ${\beta}_{k}$ is the $k$'th {\em total
Betti number} of $M$ and $p$ is the projective dimension of $M$. 

The module $M$ has a {\em pure resolution} if there are 
constants $d_0<\ldots < d_p$
such that 
$$
d_{0,i}=d_0,\ldots ,d_{p,i}=d_p 
$$
for all $i$. If in addition
$$
d_i=d_0+i,
$$
for all $1\leq i\leq p$, then 
we call the minimal free resolution  to be {\em $d_0$--linear}. 

In \cite{R} Theorem 2.7 the following bound for the Betti numbers was proved.
\begin{thm} \label{Betti_bound}
Let $M$ be an $R$--module having a pure resolution
of type $(d_0,\ldots, d_p)$ and Betti numbers $\beta_0, \ldots , \beta_p$, where 
$p$ is the projective dimension of $M$.
Then 
\begin{equation} \label{bound_Betti}
\beta_i\geq {p\choose i}
\end{equation}
for each $0\leq i\leq p$.
\end{thm}

\subsection{Hilbert--Serre Theorem}

Let $M=\bigoplus_{i\geq 0} M_i$ be a finitely generated 
nonnegatively graded module over the polynomial ring $R$.
We call the formal power series 
$$
H_M(z):= \sum_{i=0}^{\infty} h_M(i)z^i
$$
the {\em Hilbert--series} of the module $M$.

The Theorem of Hilbert--Serre states that there exists 
a (unique) polynomial $P_M(z)\in \Q[z]$, the so-called {\em Hilbert polynomial}
of $M$, such that $h_M(i)=P_M(i)$ for each $i>>0$. Moreover, $P_M$ has degree
$\mbox{dim }M-1$ and $(\mbox{dim }M-1)!$ times the leading coefficient 
of $P_M$ is the {\em Hilbert--Samuel multiplicity} of $M$, denoted here 
by $e(M)$. 

Hence there exist integers $m_0,\ldots, m_{d-1}$ 
such that $h_M(z)=m_0\cdot{z\choose d-1}+m_1\cdot{z\choose d-2}+\ldots + m_{d-1}$,
where ${z\choose r}=\frac{1}{r!}z(z-1)\ldots (z-r+1)$ and $d:=\mbox{dim} M$. 
Clearly $m_0=e(M)$.

We can summarize the Hilbert-Serre theorem as follows: 
\begin{thm} (Hilbert--Serre) \label{Hilbert_Serre}
Let $M$ be a finitely generated nonnegatively graded 
$R$--module of dimension $d$, then the following stetements hold:\\
(a) There exists a (unique) polynomial $P(z)\in \Z[z]$ such that 
the Hilbert--series $H_M(z)$ of $M$ may be written as 
$$
H_M(z)=\frac{P(z)}{(1-z)^d}
$$
(b) $d$ is the least integer for which $(1-z)^dH_M(z)$ is a polynomial.
\end{thm}


\subsection{Simplicial complexes and Stanley--Reisner rings}

We say that $\Delta\subseteq 2^{[n]}$ is a {\em simplicial complex}
 on the vertex set $[n]=\{1,2,\ldots ,n\}$, if 
 $\Delta$ is a set of subsets of $[n]$ such that 
$\Delta$ is a down--set, that is, $G\in\Delta$ and $F\subseteq G$ 
implies that $F\in \Delta$, and $\{i\}\in \Delta$ for all $i$.

The elements of $\Delta$ are called {\em faces} 
and the {\em dimension} of a face is one less than its cardinality. An $r$-face is an abbreviation for an $r$-dimensional face.
The dimension of $\Delta$ is the dimension of a maximal face.
We use the notation $\mbox{dim}(\Delta)$ for the dimension 
of $\Delta$.

If $\mbox{dim}(\Delta)=d-1$, then 
the $(d+1)$--tuple $(f_{-1}(\Delta),\ldots ,f_{d-1}(\Delta))$ 
is called the {\em $f$-vector} of $\Delta$, where
 $f_i(\Delta)$ denotes the number of $i$--dimensional faces 
of $\Delta$.

Let $\Delta$ be an arbitrary simplicial complex
on $[n]$. The {\em Stanley--Reisner ring} 
$k[\Delta]:=R/I(\Delta)$ of $\Delta$
is the quotient of the ring $R$ by the
 {\em Stanley--Reisner ideal}
$$
I(\Delta):=\langle x^F:~ F \notin \Delta \rangle, 
$$ 
generated by the non--faces of $\Delta$.

\begin{prop} \label{multi}
Let $\Delta$ be a $(d-1)$--dimensional simplicial complex. Then 
$$
e(k[\Delta])=f_{d-1}.
$$
\end{prop}
\proof 
It follows from \cite{BH} Proposition 4.1.9 and (\ref{Hilbert10}) that
$$
e(k[\Delta])=\left(\sum_{i=0}^d h_iz^i\right)\mid_{z=1}=\sum_{i=0}^d h_i=f_{d-1}.
$$
The following Theorem was proved in \cite{BH} Theorem 5.1.7. \qed

\begin{thm} \label{Hilb-Serre}
Let $\Delta$ be a $d-1$--dimensional simplicial complex with $f$-vector 
$f(\Delta):=(f_{-1},\ldots ,f_{d-1})$. 
Then the Hilbert--series 
of the Stanley--Reisner ring $k[\Delta]$ is
$$
H_{k[\Delta]}(z)=\sum_{i=-1}^{d-1} \frac{f_it^{i+1}}{(1-t)^{i+1}}.
$$
\end{thm}
\qed

Recall from Theorem \ref{Hilbert_Serre}
that a homogeneous $k$-algebra $M$ of dimension $d$ has 
a Hilbert series of the forn 
$$
H_M(z)=\frac{P(z)}{(1-z)^d}
$$
where $P(z)\in \Z[z]$. 
Let $\Delta$ be a $(d-1)$--dimensional simplicial complex and write 
\begin{equation} \label{Hilbert10}
H_{k[\Delta]}(z)=\frac{\sum_{i=0}^d h_iz^i}{(1-z)^d}.
\end{equation}

\begin{lemma} \label{f_h_vector}
The $f$-vector and the $h$-vector of a $(d-1)$--dimensional simplicial complex
$\Delta$ are related by
$$
\sum_{i} h_it^i=\sum_{i=0}^d f_{i-1}t^i(1-t)^{d-i}.
$$
In particular, the $h$-vector has length at most $d$, and 
$$
h_j=\sum_{i=0}^j (-1)^{j-i} {d-i \choose j-i}f_{i-1},\ \ f_{j-1}=\sum_{i=0}^j {d-i \choose j-i}h_i
$$
for each $j=0,\ldots ,d$ and
\end{lemma}
\qed

Let $\Gamma$ be a $(d-1)$--dimensional simplicial complex. Suppose 
that the Stanley--Reisner ring $k[\Gamma]$ has a pure free resolution
\begin{equation}
{{\cal F}}_{\Gamma}: 0\longrightarrow R(-d_p)^{\beta_p} \longrightarrow \ldots \longrightarrow 
\end{equation}

\begin{equation} \label{free99}
\longrightarrow R(-d_1)^{\beta_1} \longrightarrow R(-d_0)^{\beta_0} \longrightarrow 
R \longrightarrow k[\Gamma] \longrightarrow 0.         
\end{equation}
Here $p$ is the projective dimension of the Stanley--Reisner ring $k[\Delta]$.

It was proved in \cite{HH} Corollary 6.1.7. that
\begin{equation} \label{mult}
(-1)^{n-d}(n-d)!e(k[\Gamma])=\sum_{i=0}^p (-1)^i \beta_i (d_i)^{n-d}.
\end{equation}

The following Theorem was proved in \cite{ER} Theorem 3.
\begin{thm}(Eagon--Reiner) \label{Eagor-R}
Let $k$ be an arbitrary field and $\Delta$ be a simplicial complex on $[n]$. 
Then $I(\Delta^*)$ has a $q$--linear resolution if and only if the Stanley--Reisner ring 
$k[{\Delta}]$ is Cohen--Macaulay of dimension $n-q$..
\end{thm}


\subsection{Alexander dual}

Let $\cF\subseteq 2^{[n]}$ be an arbitrary set system. Define the complement of $\cF$ as
$$
{\cF}':= 2^{[n]}\setminus \cF.
$$
Consider the following set system
$$
\mbox{co}(\cF):=\{[n]\setminus F:~ F\in \cF\}.
$$
We denote by $\cF^*$ the {\em Alexander dual  of $\cF$} 
$$
\cF^*:=\mbox{co}({\cF}')=(\mbox{co}({\cF}))'\subseteq 2^{[n]}.
$$

Let $\Delta^*$ denote the Alexander dual of the simplicial complex $\Delta$.
We can easily compute $f(\Delta^*)$, the $f$-vector of $\Delta^*$:
\begin{lemma} \label{Alex_dual_fvector}
Let $f(\Delta)=(f_{-1}(\Delta),\ldots ,f_{d-1}(\Delta))$ be the $f$--vector of a $(d-1)$--dimensional simplicial complex $\Delta$. Then the $f$--vector of the simplicial complex $\Delta^*$ is: 
\begin{equation}
f(\Delta^*)=[\underbrace{1}_{f^*_{-1}},\underbrace{n\choose 1}_{f^*_{0}},\ldots, ,\underbrace{n\choose n-d-1}_{f^*_{n-d-1}},\underbrace{{n\choose n-d}-f_{d-1}}_{f^*_{n-d}},\ldots, \underbrace{{n\choose 2}-f_1}_{f^*_{n-2}}].
\end{equation}
Hence
$$
f_i^*={n\choose i+1}-f_{n-i-1}
$$
for each $-1\leq i\leq n-1$.
\end{lemma} 
\begin{cor} \label{dim_Alex}
Let $\Delta$ be any simplicial complex. Let $k^*$ denote the positive integer such that
${n \choose k^*}\neq f_{k^*-1}$, but ${n \choose k^*-1}=f_{k^*-2}$.
Then $d^*=\mbox{dim}(k[\Delta^*])=n-k^*$.
\end{cor}


\section{Our main result}
\subsection{Pure resolutions}


In the following Theorem we describe 
the Betti numbers of $k[\Delta]$ in terms of the $h$--vector 
of $\Delta$.
\begin{thm} \label{main_Betti}
Let $\Delta$ be a $(d-1)$--dimensional simplicial complex. Suppose 
that the Stanley--Reisner ring $k[\Delta]$ has a pure free resolution
\begin{equation}
{{\cal F}}_{\Delta}: 0\longrightarrow R(-d_p)^{\beta_p} 
\longrightarrow \ldots \longrightarrow 
\end{equation}

\begin{equation} \label{free6}
\longrightarrow R(-d_1)^{\beta_1} \longrightarrow R(-d_0)^{\beta_0} \longrightarrow 
R \longrightarrow k[\Delta] \longrightarrow 0.         
\end{equation}
Here $p$ is the projective dimension of the Stanley--Reisner ring $k[\Delta]$.

If $h(\Delta):=(h_{0}(\Delta),\ldots ,h_{d}(\Delta))$ is the 
$h$-vector of the complex $\Delta$, then
$$
\beta_i=\sum_{\ell=0}^{d_i} (-1)^{\ell+i+1}{n-d \choose \ell}h_{d_i-\ell}
$$
for each $0\leq i\leq p$.
\end{thm}
\begin{remark}
Clearly $h_i=0$ for each $i>d$.
\end{remark}

\begin{remark}
J. Herzog and M. K\"uhl proved similar formulas for the Betti number in \cite{HK} 
 Theorem 1. Here we did not assume that the Stanley--Reisner ring $k[\Delta]$
 with pure resolution is Cohen--Macaulay.   
\end{remark}
\proof
Let $M:=k[\Delta]$ denote the Stanley--Reisner ring of $\Delta$. 
Then we infer from Theorem \ref{Hilb-Serre} that
\begin{equation} \label{Hilbert8}
H_M(z)=\frac{\sum_{i=0}^d h_iz^i}{(1-z)^d}.
\end{equation}
Since the Hilbert--series is additive on short exact sequences,
and since 
$$
H_R(z)=\frac{1}{(1-z)^n},
$$
and consequently
$$
H_{R(-s)}(z)=\frac{z^s}{(1-z)^n},
$$
the pure resolution
\begin{equation}
{{\cal F}}_{\Delta}: 0\longrightarrow R(-d_p)^{\beta_p} \longrightarrow \ldots \longrightarrow 
\end{equation}

\begin{equation} \label{free5}
\longrightarrow R(-d_1)^{\beta_1} \longrightarrow 
R(-d_0)^{\beta_0} \longrightarrow R \longrightarrow M \longrightarrow 0.         
\end{equation}
yields to
\begin{equation} \label{Hilbert7}
H_M(z)=\frac{1}{(1-z)^n}+\sum_{i=0}^p (-1)^{i+1} \beta_i \frac{z^{d_i}}{(1-z)^n}, 
\end{equation}
where $p=pdim(M)$.

Write $d:=\mbox{dim} M$, and let $m:=\mbox{codim}(M)=n-d$. 
It follows from the Auslander--Buchbaum formula that $m\leq p$. 

Comparing the two expressions (\ref{Hilbert7}) and (\ref{Hilbert8})
for $H_M$, we find
\begin{equation} \label{Hilbert6}
(1-z)^m \left(\sum_{i=0}^d h_iz^i \right)=\sum_{i=0}^p (-1)^{i+1} \beta_i z^{d_i}+1
\end{equation}

Using the binomial Theorem we get that
\begin{equation} \label{Hilbert9}
\left( \sum_{j=0}^{n-d} (-1)^j{n-d\choose j}z^j \right)\left(\sum_{i=0}^d h_iz^i \right)
=\sum_{i=0}^p (-1)^{i+1} \beta_i z^{d_i}+1
\end{equation}

Comparing the coefficients on the two sides of (\ref{Hilbert9}),
 we get the result.
\qed
\begin{cor} \label{main_cor_h_vect}
Let $\Delta$ be a $(d-1)$--dimensional simplicial complex. Suppose 
that the Stanley--Reisner ring $k[\Delta]$ has a pure free resolution
\begin{equation}
{{\cal F}}_{\Delta}: 0\longrightarrow R(-d_p)^{\beta_p} \longrightarrow \ldots \longrightarrow 
\end{equation}

\begin{equation} \label{free66}
\longrightarrow R(-d_1)^{\beta_1} \longrightarrow R(-d_0)^{\beta_0} \longrightarrow 
R \longrightarrow k[\Delta] \longrightarrow 0.         
\end{equation}
Here $p$ is the projective dimension of the Stanley--Reisner ring $k[\Delta]$.
Let $h(\Delta):=(h_{0}(\Delta),\ldots ,h_{d}(\Delta))$ denote the 
$h$-vector of the complex $\Delta$. Let $s$ be a positive integer such that $s\neq d_i$
for each $0\leq i\leq p$. Then
$$
\sum_{\ell=0}^{s} (-1)^{\ell}h_{s-\ell}{n-d \choose \ell}=0.
$$
\end{cor}
\proof
Define
$$
P(z):=1+\sum_{i=0}^p (-1)^{i+1} \beta_i z^{t+i}\in {\Q}[z].
$$
Comparing the
coefficients of both side of (\ref{Hilbert9}), we get the result. \qed

\begin{cor} \label{main_cor2}
Let $\Delta$ be a $(d-1)$--dimensional simplicial complex. Suppose 
that the Stanley--Reisner ring $k[\Delta]$ has a pure free resolution
\begin{equation}
{{\cal F}}_{\Delta}: 0\longrightarrow R(-d_p)^{\beta_p} \longrightarrow 
\ldots \longrightarrow 
\end{equation}

\begin{equation} \label{free9}
\longrightarrow R(-d_1)^{\beta_1} \longrightarrow R(-d_0)^{\beta_0} 
\longrightarrow R \longrightarrow k[\Delta] \longrightarrow 0.         
\end{equation}
Here $p$ is the projective dimension of the Stanley--Reisner ring $k[\Delta]$.
Then
\begin{equation}
\sum_{\ell=0}^{d_i} (-1)^{\ell+i+1}{n-d \choose \ell}h_{d_i-\ell}\geq {p\choose i}
\end{equation}
for each $0\leq i\leq p$.
\end{cor}
\proof 
This follows easily from Theorem \ref{Betti_bound} and Theorem \ref{main_Betti}.
\qed

\subsection{$t$--linear free resolutions}

In the following we specialize our results to $t$--linear resolutions. 
Since any $t$--linear free resolution will be a pure resolution, so we can apply
Theorem \ref{main_Betti}.

\begin{cor} \label{main_cor}
Let $\Delta$ be a $(d-1)$--dimensional simplicial complex. Suppose 
that the Stanley--Reisner ring $k[\Delta]$ has a $t$--linear free resolution
\begin{equation}
{{\cal F}}_{\Delta}: 0\longrightarrow R(-t-p)^{\beta_p} \longrightarrow \ldots \longrightarrow 
\end{equation}
\begin{equation} \label{free7}
\longrightarrow R(-t-1)^{\beta_1} \longrightarrow R(-t)^{\beta_0} 
\longrightarrow R \longrightarrow k[\Delta] \longrightarrow 0.         
\end{equation}
Here $p$ is the projective dimension of the Stanley--Reisner ring $k[\Delta]$.
If $h(\Delta):=(h_{0}(\Delta),\ldots ,h_{d}(\Delta))$ is the 
$h$-vector of the complex $\Delta$, then
$$
\beta_i=\sum_{\ell=0}^{t+i} (-1)^{\ell+i+1}h_{t+i-\ell}{n-d \choose \ell}
$$
for each $0\leq i\leq p$.
\end{cor}

\begin{cor} \label{main_cor3}
Let $\Delta$ be a $(d-1)$--dimensional simplicial complex. Suppose 
that the Stanley--Reisner ring $k[\Delta]$ has a $t$--linear free resolution
\begin{equation}
{{\cal F}}_{\Delta}: 0\longrightarrow R(-t-p)^{\beta_p} \longrightarrow \ldots \longrightarrow 
\end{equation}

\begin{equation} \label{free14}
\longrightarrow R(-t-1)^{\beta_1} \longrightarrow R(-t)^{\beta_0} 
\longrightarrow R \longrightarrow k[\Delta] \longrightarrow 0.         
\end{equation}
Here $p$ is the projective dimension of the Stanley--Reisner ring $k[\Delta]$.

If $h(\Delta):=(h_{0}(\Delta),\ldots ,h_{d}(\Delta))$ is the 
$h$-vector of the complex $\Delta$, then
$$
\sum_{\ell=0}^{j} (-1)^{\ell}h_{j-\ell}{n-d \choose \ell}=0.
$$
for each $p+t<j\leq n$ and $0<j<t$.
\end{cor}
\begin{cor} \label{main_cor_inequ}
Let $\Delta$ be a $(d-1)$--dimensional simplicial complex. Suppose 
that the Stanley--Reisner ring $k[\Delta]$ has a $t$--linear free resolution
\begin{equation}
{{\cal F}}_{\Delta}: 0\longrightarrow R(-t-p)^{\beta_p} \longrightarrow \ldots \longrightarrow 
\end{equation}

\begin{equation} \label{free12}
\longrightarrow R(-t-1)^{\beta_1} \longrightarrow R(-t)^{\beta_0} 
\longrightarrow R \longrightarrow k[\Delta] \longrightarrow 0.         
\end{equation}
Here $p$ is the projective dimension of the Stanley--Reisner ring $k[\Delta]$.
If $h(\Delta):=(h_{0}(\Delta),\ldots ,h_{d}(\Delta))$ is the 
$h$-vector of the complex $\Delta$, then
\begin{equation}
\sum_{\ell=0}^{t+i} (-1)^{\ell+i+1}{n-d \choose \ell}h_{t+i-\ell}\geq {p\choose i}
\end{equation}
for each $0\leq i\leq p$.
\end{cor}

\section{Applications}
\subsection{Cohen--Macaulay squarefree monomial ideals}
Now we apply our main results for Cohen--Macaulay simplicial complexes.
\begin{cor} \label{Cohen--Macaulay} 
Let $\Delta$ be a $(d-1)$--dimensional Cohen--Macaulay simplicial complex.
Consider the $t$--linear free resolution
\begin{equation}
{{\cal F}}_{\Delta^*}: 0\longrightarrow R(-t-p)^{\beta^*_p} \longrightarrow \ldots \longrightarrow 
\end{equation}
\begin{equation} \label{Alex3}
\longrightarrow R(-t-1)^{\beta^*_1} \longrightarrow R(-t)^{\beta^*_0} 
\longrightarrow R \longrightarrow k[\Delta^*] \longrightarrow 0.         
\end{equation}
Here $t=n-d$. Let $k^*$ denote the positive integer such that
${n \choose k^*}\neq f_{k^*-1}$, but ${n \choose k^*-1}=f_{k^*-2}$. Then 
\begin{equation}
\sum_{\ell=0}^{j} (-1)^{\ell}h^*_{j-\ell}{k^* \choose \ell}=0.
\end{equation}
for each $p+t<j\leq n$ and $0<j<t$.
\end{cor}
\proof 
This follows easily from Corollary \ref{dim_Alex} and Corollary \ref{main_cor3}.
\qed
\begin{cor} \label{Cohen--Macaulay_Betti}
Let $\Delta$ be a $(d-1)$--dimensional Cohen--Macaulay simplicial complex. 
Consider the $t$--linear free resolution
\begin{equation}
{{\cal F}}_{\Delta^*}: 0\longrightarrow R(-t-p)^{\beta^*_p} \longrightarrow \ldots \longrightarrow 
\end{equation}
\begin{equation} \label{Alex2}
\longrightarrow R(-t-1)^{\beta^*_1} \longrightarrow R(-t)^{\beta^*_0} 
\longrightarrow R \longrightarrow k[\Delta^*] \longrightarrow 0.         
\end{equation}
Here $t=n-d$. Let $p$ denote the projective dimension of $k[\Delta^*]$. 
Let $k^*$ denote the positive integer such that
${n \choose k^*}\neq f_{k^*-1}$, but ${n \choose k^*-1}=f_{k^*-2}$.
Let $h(\Delta^*):=(h_{0}^*(\Delta),\ldots ,h_{n-2}^*(\Delta))$ be the 
$h$-vector of the complex $\Delta^*$,
Then 
\begin{equation} \label{dual_Betti}
\beta_i^*=\sum_{\ell=0}^{t+i} (-1)^{\ell+i+1}h^*_{t+i-\ell}{k^* \choose \ell}
\end{equation}
for each $0\leq i\leq p$.
\end{cor}
\begin{remark}
Eagon and Reiner proved similar formulas for the Betti numbers $\beta^*_i$
in \cite{ER} Theorem 4. They described the folowing equation:
$$
\sum_{i\geq 1} \beta_i^* t^{i-1}=\sum_{i=0}^d h_i(\Delta)(t+1)^i
$$
\end{remark}
\proof
It follows from the Eagon--Reiner Theorem \ref{Eagor-R} that
the free resolution (\ref{Alex2}) exists. Now we can apply 
Corollary \ref{main_cor} for the $t$--linear 
free resolution (\ref{Alex2}). Since $d^*=n-k^*$ by Corollary \ref{dim_Alex}, 
hence the result follows. \qed
\begin{cor} \label{Cohen--Macaulay_mult}
Let $\Delta$ be a $(d-1)$--dimensional Cohen--Macaulay simplicial complex.
Consider the $t$--linear free resolution
\begin{equation}
{{\cal F}}_{\Delta^*}: 0\longrightarrow R(-t-p)^{\beta^*_p} \longrightarrow \ldots \longrightarrow 
\end{equation}
\begin{equation} \label{Alex5}
\longrightarrow R(-t-1)^{\beta^*_1} \longrightarrow R(-t)^{\beta^*_0} 
\longrightarrow R \longrightarrow k[\Delta^*] \longrightarrow 0.         
\end{equation}
Here $t=n-d$. 
Let $k^*$ denote the positive integer such that
${n \choose k^*}\neq f_{k^*-1}$, but ${n \choose k^*-1}=f_{k^*-2}$.
Then 
\begin{equation} \label{dual_mult2}
e(k[\Delta^*])=\frac{(-1)^{k^*}}{(k^*)!}\sum_{i=0}^p (-1)^i \beta_i^* (n-d+i)^{k^*}
\end{equation}
where $p$ is the projective dimension of $k[\Delta^*]$.
\end{cor}
\proof
Let $\Delta$ be a $(d-1)$--dimensional simplicial complex.
Let $\Gamma:=\Delta^*$. Then using Eagon--Reiner Theorem \ref{Eagor-R} and 
(\ref{mult}), we get our result. \qed 

\begin{cor} \label{Cohen--Macaulay_equation}
Let $\Delta$ be a $(d-1)$--dimensional Cohen--Macaulay simplicial complex.
Consider the $t$--linear free resolution
\begin{equation}
{{\cal F}}_{\Delta^*}: 0\longrightarrow R(-t-p)^{\beta^*_p} \longrightarrow \ldots \longrightarrow 
\end{equation}

\begin{equation} \label{Alex56}
\longrightarrow R(-t-1)^{\beta^*_1} \longrightarrow R(-t)^{\beta^*_0} 
\longrightarrow R \longrightarrow k[\Delta^*] \longrightarrow 0.         
\end{equation}
Here $t=n-d$. Let $k^*$ denote the positive integer such that
${n \choose k^*}\neq f_{k^*-1}$, but ${n \choose k^*-1}=f_{k^*-2}$.
Then
\begin{equation} \label{dual_mult3}
{n\choose k^*}-f_{k^*-1}=\frac{(-1)^{k^*}}{(k^*)!}\sum_{i=0}^p (-1)^i \beta_i^* (n-d+i)^{k^*}.
\end{equation}
\end{cor}
\proof 
Since by Lemma \ref{multi}
\begin{equation} 
e(k[\Delta^*])=f^*_{d^*-1}={n\choose k^*}-f_{k^*-1}
\end{equation}
and we infer from Corollary \ref{Cohen--Macaulay_mult} that
\begin{equation} 
e(k[\Delta^*])=\frac{(-1)^{k^*}}{(k^*)!}\sum_{i=0}^p (-1)^i \beta_i^* (t+i)^{k^*}
=\frac{(-1)^{k^*}}{(k^*)!}\sum_{i=0}^p (-1)^i \beta_i^* (n-d+i)^{k^*},
\end{equation}
we get the result.
\qed

\begin{cor} \label{dual_cor_inequ}
Let $\Delta$ be a $(d-1)$--dimensional Cohen--Macaulay simplicial complex.
Consider the $t$--linear free resolution
\begin{equation}
{{\cal F}}_{\Delta^*}: 0\longrightarrow R(-t-p)^{\beta^*_p} \longrightarrow \ldots \longrightarrow 
\end{equation}

\begin{equation} \label{Alex55}
\longrightarrow R(-t-1)^{\beta^*_1} \longrightarrow R(-t)^{\beta^*_0} 
\longrightarrow R \longrightarrow k[\Delta^*] \longrightarrow 0.         
\end{equation}
Here $t=n-d$. Let $k^*$ denote the positive integer such that
${n \choose k^*}\neq f_{k^*-1}$, but ${n \choose k^*-1}=f_{k^*-2}$.
If $h(\Delta^*):=(h_{0}(\Delta^*),\ldots ,h_{d}(\Delta^*))$ is the 
$h$-vector of the complex $\Delta^*$, then
\begin{equation}
\sum_{\ell=0}^{t+i} (-1)^{\ell+i+1}{k^* \choose \ell}h^*_{t+i-\ell}\geq {p\choose i}
\end{equation}
for each $0\leq i\leq p$.
\end{cor}
\proof 
This follows easily from Theorem \ref{Eagor-R} and Corollary \ref{main_cor_inequ}.
\qed

\subsection{Chordal graphs}
We can specialize our results for chordal graphs.
\begin{cor} \label{main_gr} 
Let $G$ be an arbitrary chordal graph. Let $\Delta:=\Delta(G)$ denote the clique
 complex of
$G$ and $d:=\mbox{dim}(\Delta)+1$. 
Let $h(\Delta):=(h_{0}(\Delta),\ldots ,h_{d}(\Delta))$ denote the 
$h$-vector of the complex $\Delta$. Let $p$ be the projective 
dimension of the Stanley--Reisner ring $k[\Delta]$. Then
$$
\sum_{\ell=0}^{j} (-1)^{\ell}h_{j-\ell}{n-d \choose \ell}=0
$$
if either $j=1$ or $p+2<j\leq n$.
\end{cor}
\proof 
Let $H$ denote the complementary graph of $G$.
Then the ideal $I(H)$ has a $2$--linear free resolution 
by Theorem \ref{Froberg_theorem}. So we can apply
Corollary \ref{main_cor3}.
\qed
\begin{cor} \label{main_gr_Betti} 
Let $G$ be an arbitrary chordal graph. Let $\Delta:=\Delta(G)$ denote the clique
 complex of
$G$ and $d:=\mbox{dim}(\Delta)+1$. 
Then the Stanley--Reisner ring $k[\Delta]$ has an $2$--linear free resolution
\begin{equation}
{{\cal F}}_{\Delta}: 0\longrightarrow R(-2-p)^{\beta_p} \longrightarrow \ldots \longrightarrow 
\end{equation}

\begin{equation} \label{free77}
\longrightarrow R(-3)^{\beta_1} \longrightarrow R(-2)^{\beta_0} 
\longrightarrow R \longrightarrow k[\Delta] \longrightarrow 0.         
\end{equation}
Here $p$ is the projective dimension of the Stanley--Reisner ring $k[\Delta]$.
If $h(\Delta):=(h_{0}(\Delta),\ldots ,h_{d}(\Delta))$ is the 
$h$-vector of the complex $\Delta(G)$, then
$$
\beta_i=\sum_{\ell=0}^{2+i} (-1)^{\ell+i+1}h_{2+i-\ell}{n-d \choose \ell}
$$
for each $0\leq i\leq p$.
\end{cor}
\proof
Let $H$ denote the complementary graph of $G$.
Then the ideal $I(H)$ has a $2$--linear free resolution 
by Theorem \ref{Froberg_theorem}. Corollary \ref{main_cor}
gives the result.

\begin{cor} \label{chordal_cor_inequ}
Let $G$ be an arbitrary chordal graph. Let $\Delta:=\Delta(G)$ denote the clique
 complex of
$G$ and $d:=\mbox{dim}(\Delta)+1$. 
Let $p$ be the projective dimension of the Stanley--Reisner ring $k[\Delta]$.
If $h(\Delta):=(h_{0}(\Delta),\ldots ,h_{d}(\Delta))$ is the 
$h$-vector of the complex $\Delta$, then
\begin{equation}
\sum_{\ell=0}^{2+i} (-1)^{\ell+i+1}{n-d \choose \ell}h_{2+i-\ell}\geq {p\choose i}
\end{equation}
for each $0\leq i\leq p$.
\end{cor}
\proof
Let $H$ denote the complementary graph of $G$.
Then the ideal $I(H)$ has a $2$--linear free resolution 
by Theorem \ref{Froberg_theorem}. The result follows from
 Corollary \ref{main_cor_inequ}

{\bf Acknowledgements.}  I am indebted to Josef Schicho, Russ Woodroofe 
 and Lajos R\'onyai 
for their useful remarks.

\end{document}